\documentclass[12pt]{article}
\usepackage[mathscr]{eucal}
\usepackage{amssymb,latexsym}
\usepackage{verbatim}
\usepackage{amsmath}
\usepackage{amsthm}
\usepackage{enumerate}

\newtheorem{thm}{Theorem}[section]

\newtheorem{prop}[thm]{Proposition}
\newtheorem{cor}[thm]{Corollary}

\newcommand\calA{{\mathcal{A}}}
\newcommand\calB{{\mathcal{B}}}

\newcommand\calV{{\mathcal{V}}}

\newcommand\calT{{\mathcal{T}}}

\newcommand\calH{{\mathcal{H}}}

\renewcommand\l{\lambda}

\newcommand\bbH{{\mathbb{H}}}
\newcommand\bbR{{\mathbb R}}

\renewcommand\S{\Sigma}

\renewcommand\d{\partial}

\newcommand\f{\phi}
\renewcommand\L{\triangle}
\newcommand\D{\nabla}
\newcommand\e{\epsilon}

\renewcommand\div{{\rm div}}

\newcommand\la{\langle}
\newcommand\ra{\rangle}

\renewcommand\l{\lambda}
\newcommand\g{\gamma}

\renewcommand\th{\theta}

\newcommand\<{\la}
\renewcommand\>{\ra}

\newcommand\beq{\begin{equation}}
\newcommand\eeq{\end{equation}}
\newcommand\ben{\begin{enumerate}}
\newcommand\een{\end{enumerate}}
\newcommand\bit{\begin{itemize}}
\newcommand\eit{\end{itemize}}

\newcounter{mnotecount}[section]

\title{Stability and rigidity of extremal surfaces in Riemannian geometry and General Relativity}

\author{Gregory J. Galloway\thanks{Supported in part by NSF grant DMS-0708048.}
 \\Department of Mathematics \\ University of Miami
}

\begin{document}
\date{}
\maketitle
\vspace{.2in}

\begin{abstract} 
We examine some common features of minimal surfaces, nonzero constant mean curvature surfaces and marginally outer trapped surfaces, concerning their stability and rigidity, and consider some applications to Riemannian geometry and general relativity.

\footnotesize{


\vspace*{4mm}
\noindent{\bf 2000 Mathematics Subject Classification: 53C20, 83C57}

\vspace*{2mm}
\noindent{\bf Keywords and Phrases: minimal surfaces, marginally trapped surfaces, black hole topology, positive mass theorems} 
}

\end{abstract}


\section{Minimal hypersurfaces in manifolds of nonnegative scalar curvature}

It is an honor to contribute a paper to this volume celebrating Rick Schoen's 60th birthday.
Rick's work has had a tremendous  
influence on the research of so many of us
working in the areas of differential geometry, geometric analysis and mathematical relativity.
It is a pleasure to consider in this article  some work of the author and others in which the impact
of Rick's work is so strongly  evident.  

In this first section we review some fundamental results of Schoen and Yau, along with
some refinements of those results, concerning compact minimal hypersurfaces in manifolds of 
nonnegative scalar curvature.  In the next section we discuss joint work with Schoen, and some 
refinements of that work,  which show how some of the results described in this section extend to marginally
outer trapped surfaces in general initial data sets satisfying the dominant energy condition.  The results in this context yield a generalization of Hawking's black hole topology theorem to higher dimensions.  In
Section 3, we describe  a positive mass result for asymptotically hyperbolic manifolds, obtained in joint work with Andersson and  Cai, which is proved via the general ``minimal surface methodology" of Schoen and Yau,  and which, hence, does not require the assumption of spin.  The proof makes use of natural extensions of the results discussed in this section to positive constant mean curvature hypersurfaces in manifolds with negative lower bound on the  scalar curvature.

In the seminal papers \cite{SYsc, SYsc2}, Schoen and Yau studied the question of
which compact manifolds admit a metric of positive scalar curvature.  Prior to this work,
Lichnerowicz  \cite{Lich} had proved, using a Bochner type argument associated with the Dirac operator
and the Atiyah-Singer index theorem,  that compact $4k$-dimensional spin manifolds of positive scalar curvature must have vanishing $\hat A$-genus.  This was followed later by work of Hitchin \cite{Hitchin}, who also used the  spinorial method to obtain further obstructions to the existence of metrics of positive scalar curvature.  While these results are quite striking, they left open some very basic questions, for example, the question as to whether the $k$-torus, $k \ge 3$, admits a metric of positive scalar curvature. 
Then in \cite{SYsc}, Schoen and Yau made a huge advance by proving, using minimal surface techniques, that if the fundamental group of a compact orientable $3$-manifold  contains a subgroup isomorphic to a nontrivial surface group
then the manifold does not admit a metric of positive scalar curvature. This implies, of course, that
the $3$-torus does not admit a metric of positive scalar curvature.   In \cite{SYsc2}, Schoen and Yau generalized their techniques to higher dimensions, thereby establishing inductively the existence of a large class of compact manifolds, including tori, of dimension up to~$7$, that do not admit metrics of positive scalar curvature.   The fundamental obeservation made in   \cite{SYsc2} is that if $M^n$,  $3 \le n \le 7$, is a compact orientable manifold of positive scalar curvature then any nontrivial 
codimension-one homology class can be represented by a manifold that admits a metric of positive scalar curvature.  This is proved by choosing a manifold of least area in the homology class, and making use of the positivity of the stability operator, ``rearranged"  in an especially useful way.


We would like to  review this argument in a bit more detail.  Let $\S^{n-1}$ be a compact two-sided
hypersurface in an $n$-dimensional Riemannian manifold $(M^n,g)$.  With respect to a chosen
unit normal field $\nu$ along $\S$,  introduce the second fundamental form $B$ of  $\S$,
$B(X,Y) = \<\D_X\nu,Y\>$,  $X,Y \in T_p\S$, $ p\in \S$, and  mean curvature $H = {\rm tr}_hB$,
where $h$ is the induced metric on $\S$.

Now assume $\S$ is minimal, $H = 0$, and consider normal variations of $\S$, i.e.,  variations   
$t \to \S_t$ of $\S = \S_0$, $-\e <t < \e$, with variation vector field 
$V = \left . \frac{\d}{\d t}\right |_{t=0} = \phi N$,  $\phi \in C^{\infty}(\S)$.  Let $A(t) =$ area of $\S_t$; then
by the formula for the second variation of area,
\beq
A''(0) =  \int_{\S} \phi L(\phi) dA  \,,
\eeq 
where $L : C^{\infty}(\S) \to C^{\infty}(\S)$ is the operator,
\beq\label{stabop0}
L(\phi) = -\triangle \phi - ({\rm Ric}(\nu,\nu) + |B|^2) \phi \,.
\eeq

By definition, $\S$ is stable provided $A''(0) \ge 0$, for all normal variations of $\S$.  
From Rayleigh's formula for the principal eigenvalue $\l_1(L)$ of the {\it stability operator} $L$,
we obtain the following well-known fact.
\begin{prop}\label{stable}
The following three conditions are equivalent: (i) $\S$ is stable.  (ii)~$\lambda_1(L) \ge 0$. (iii) There exists $\f \in C^{\infty}(\S)$, $\f > 0$, such that $L(\phi) \ge~0$.
\end{prop}

\medskip
The key geometric observation made by Schoen and Yau in \cite{SYsc2} is the following.

\begin{prop}\label{pos} 
Let $(M^n,g)$, $n \ge 3$ be a Riemannian manifold with nonnegative scalar curvature,
$S \ge 0$.  If $\S^{n-1}$ is a compact two-sided stable minimal hypersurface in $M^n$ then $\S^{n-1}$ admits a metric of positive scalar curvature, unless it is totally geodesic and Ricci flat, and $S = 0$
along $\S^{n-1}$.
\end{prop}

\proof  We give here a slight variation of the proof in \cite{SYsc2}.  
Rather than compare the stability operator to the conformal Laplacian (which then requires
$n \ge 4$), it was observed in \cite{Cai}, that, by using a ``nonstandard" conformal factor, one can work directly with the stability operator.

The essential step taken by Schoen and Yau was to rewrite the stability operator in an especially useful way.   From  the Gauss equation one obtains,
\begin{align}
{\rm Ric}(\nu,\nu) + |B|^2 &= \frac12\left(S + |B|^2 + H^2 - S_{\S}  \right) \label{gauss}  \\
&=  \frac12\left(S + |B|^2  - S_{\S}  \right)
\end{align}
where $S_{\S}$ is the scalar curvature of $\S$ in the induced metric, and where in the second line
we have used that $\S$ is minimal.   Hence the stability operator $L$ may be written as,
\beq\label{stabop}
L(\phi) = -\triangle \phi + \frac12(S_{\S} - S - |B|^2) \phi \,.
\eeq

Now let $\phi$ be a positive eigenfunction associated to $\l_1(L)$, and consider $\S$ in the conformally related metric $\tilde h = \f^{\frac2{n-2}} h$.  The scalar curvature $\tilde S$ of $\S$ in the metric $\tilde h$ is given by,
\begin{align}\label{eq:conf}
\tilde S & =  \phi^{-\frac{n}{n-2}}(-2\L\f +  S_{\S} \f +\frac{n-1}{n-2}\frac{|\D\f|^2}{\f}) \nonumber \\  
& =  \phi^{-\frac{2}{n-2}}(2\l_1(L) + S + |B|^2 +\frac{n-1}{n-2}\frac{|\D\f|^2}{\f^2}) \,.
\end{align}
By our assumptions, all terms on the right hand side of the above are nonnegative, and hence
$\tilde S \ge 0$.  Then by well-known results  \cite{KW}, either $\S$ admits a metric of positive scalar curvature or  $\tilde S$ vanishes identically.  In the latter case, all terms on the right hand side of \eqref{eq:conf} vanish.  Using this in \eqref{stabop} implies that $S_{\S} \equiv 0$.  Then, by a result of Bourguignon (see \cite{KW}), $\S$ either carries a metric of positive scalar curvature or is Ricci flat.
\qed

\smallskip
The exceptional case in Proposition \ref{pos} can occur.   For example, it is not hard to construct a metric of nonnegative scalar curvature on the $3$-sphere $S^3$, such that $S^3$ in this metric contains a stable minimal torus.  The papers \cite{CG, Cai} considered the rigidity 
of Proposition \ref{pos} in the case that $\S$ is area minimizing (rather than just stable) and is {\it not} of  positive Yamabe type, i.e., does not admit a metric of positive scalar curvature.  (The issue of
such rigidity in three dimensions was raised in \cite{FCS}.) 
Since variations on the rigidity result obtained in  \cite{Cai} have been useful for certain problems in general relativity (see \cite{ACG, G}, as well as the following sections), 
we take some time here to recall this result. 
\begin{thm}\label{rigid}
Let $(M^n,g)$, $n \ge 3$ be a Riemannian manifold with nonnegative scalar curvature,
$S \ge 0$, and suppose $\S^{n-1}$ is a compact two-sided hypersurface in $M^n$ which 
locally minimizes area.  If $\S$  is not of positive Yamabe type then a neighborhood $U$ of $\S$ splits, i.e., 
$(U,g|_U)$ is isometric to $((-\e,\e)\times \S, dt^2 + h)$, where $h$, the induced metric on 
$\S$, is Ricci flat.  
\end{thm}  
 
The assumption that $\S$   ``locally minimizes"  area can be taken to mean, for example, that $\S$ has area less than or equal to that of all graphs over $\S$ with respect to a fixed Gaussian normal
coordinate system. 

\proof The original proof given in \cite{Cai} involves a deformation of the
ambient space metric  in a neighborhood of $\S$ to a metric of strictly positive scalar curvature.  
The proof  given here is patterned after the proof of Theorem 2.3 in \cite{ACG} (see also \cite{G, bray}), and does not require such a deformation.   

Let $\calH(u)$ denote the mean curvature of the hypersurface $\S_u : x \to exp_x u(x)\nu$, $u\in C^{\infty}(\S)$, $u$ sufficiently small.  $\calH$ has linearization $\calH'(0) = L$,
where $L$ is the stability operator~(\ref{stabop}). But by Proposition \ref{pos},
$L$ reduces to $-\triangle$, and hence $\calH'(0) = - \triangle$.  Since the kernel of 
$\calH'(0)$ consists only of the constants, the implicit function theorem implies that a neighborhood
of $\S$ is foliated by constant mean curvature hypersurfaces.   More explicitly, 
there exists a neighborhood $U$ of $\S$ such that, up to isometry,
\beq\label{coords}
U = (-\e, \e) \times \S \qquad g|_U = \phi^2 dt^2 + h_t \,,
\eeq
where $h_t = h_{ij}(t,x)dx^idx^j$ ($x^i$ local coordinates on $\S$), $\phi = \phi(t,x)$ and $\S_t = \{t\} \times \S$  has constant mean curvature.  Let $A(t) =$ the area of $\S_t$; since $\S$ is locally
of least area, we have that $A(0) \le A(t)$ for all $t \in (-\e,\e)$, for $\e$  small enough.

Let $H(t)$ denote the
mean curvature of $\S_t$.    $H = H(t)$ obeys the evolution
equation,
\beq\label{evolve}
\frac{dH}{dt} = - \triangle \phi + \frac12(S_{\S_t} - S - |B|^2 -H^2) \phi \,,
\eeq
(compare \eqref{stabop}).  Since $\S$ is minimal, $H(0) = 0$.  We claim that $H \le 0$ for $t \in [0,\e)$.  If not, there exists $t_0 \in (0,\e)$ such that $H'(t_0) > 0$.  Let $\tilde S$ be the scalar curvature of $\S_{t_0}$ in the conformally related metric $\tilde h = \f^{\frac2{n-2}}h_{t_0}$.  Arguing similarly as in the derivation of (\ref{eq:conf}), Equation  (\ref{evolve}) implies,
\beq\label{rescale2}
\tilde S =  \phi^{-\frac{2}{n-2}}(2\f^{-1}H'(t_0)  + S + |B|^2 + H^2 +\frac{n-1}{n-2}\frac{|\D\f|^2}{\f^2})  > 0 \,,
\eeq
contrary to the assumption that $\S_{t_0} \approx  \S$ is not of positive Yamabe type.  

Thus, $H \le 0$ on $[0,\e)$ as claimed, and hence,
\beq\label{firstvar}
A'(t)  =  \int_{\S_t}  H  \f\, dA \, \le 0 \,, \quad \mbox{for all  } t \in [0,\e) \,.
\eeq
But since $A = A(t)$ achieves a minimum at $t=0$, it must be that
$A'(t) = 0$ for $t \in [0,\e)$.  Hence,  the integral in (\ref{firstvar})
vanishes, which implies that $H = 0$ on  $[0,\e)$.  A similar argument shows that 
$H = 0$ on $(-\e, 0]$,
as well.  
Equation (\ref{evolve}) then implies that $L(\f) = 0$ on each $\S_t$, where $L$ is as in
Equation \eqref{stabop}.  It then follows  from Proposition \ref{stable}, that each $\S_t$ is stable.  
From Proposition~\ref{pos}  we have that each $\S_t$ is totally geodesic, so that  
$\frac{\d h_{ij}}{\d t} = 0$ for each $t$, and that $\f$ depends only on $t$.  By a simple change of $t$-coordinate
in (\ref{coords}), we may assume without loss of generality that $\f = 1$.
The result follows.~\qed

\proof[Remarks]
If, in the setting of Theorem \ref{rigid},  $\S$ is assumed to be minimal and locally of least area only to one side of 
$\S$  (the side into which $\nu$ points, say), then it is easy to check that one still gets a splitting locally to that side.

Now, suppose $\S$ is a compact minimal hypersurface in $M$ that separates $M$ into an ``inside" and an ``outside".  We say that $\S$ is an outermost minimal hypersurface if there are no minimal
hypersurfaces outside of and homologous to $\S$.  We say, further, that $\S$ is outer
area minimizing  if its area is less than or equal to the area of any other hypersurface outside of and homologous to it.  Given an outermost minimal hypersurface $\S$  in $M$, if there is a mean convex barrier  $\S'$ outside of and homologous to $\S$, and if dim $M \le 7$, then $\S$ will be outer area minimizing in the region $W$ bounded by $\S$ and $\S'$, as can be seen by minimizing area
in the homology class of $\S$ in $W$.  See e.g., \cite{HI, Eich2} for natural conditions which
guarantee the existence of outermost, outer area minimizing minimal hypersurfaces in Riemannian
manifolds.

Theorem \ref{rigid} and the above remarks imply the following.

\begin{cor}\label{outermost}
Let $(M^n,g)$, $n \ge 3$ be a Riemannian manifold with nonnegative scalar curvature,
$S \ge 0$.  If $\S^{n-1}$ is an outermost, outer area minimizing minimal hypersurface in $M^n$,
then $\S^{n-1}$ must be of positive Yamabe type.  
\end{cor}

\proof By Theorem \ref{rigid} and the remarks above, if $\S$ is not of positive Yamabe type then
an outer neighborhood of $\S$ will be foliated by minimal (in fact, totally geodesic) hypersurfaces
isotopic to $\S$, thereby contradicting the assumption that $\S$ is outermost.~\qed

We note that, in Corollary \ref{outermost}, if $\S$ is outer area minimizing with respect to
an asymptotically flat end, then the assumption of being outermost can be dropped.  For, in this case,
one can use Theorem \ref{rigid} to show that the entire end splits, contrary to it being
asymptotically flat. 

\smallskip
To conclude this section, we wish to mention a recent paper of  Bray, Brendle and Neves \cite{bray}, in which they obtain a rigidity result for least area $2$-spheres, whose proof is similar in spirit to that of Theorem \ref{rigid}.
Suppose  $(M,g)$ is a compact $3$-manifold with $\pi_2(M) \neq 0$. Denote by $\mathscr{F}$ the set of all smooth maps $f: S^2 \to M$ which represent a non-trivial element of $\pi_2(M)$, and define
\begin{equation} 
\label{definition.of.A}
\mathscr{A}(M,g) = \inf \{\text{\rm area}(S^2,f^* g): f \in \mathscr{F}\}. 
\end{equation}
The main result in \cite{bray} is then as follows.

\begin{thm}
\label{bbnthm}
We have
\begin{equation}
\label{upper.bound.for.area}
\mathscr{A}(M,g) \, \inf_M S \leq 8\pi, 
\end{equation}
where $S$ denotes the scalar curvature of $(M,g)$. Moreover, if the equality holds, then the universal cover of $(M,g)$ is isometric to the standard cylinder $S^2 \times \mathbb{R}$ up to scaling.
\end{thm}

Thus, rigidity can be achieved even in the positive Yamabe case.   The proof 
makes essential use of the Gauss-Bonnet theorem, and hence the restriction to three dimensions.
The authors first establish a local splitting about an area minimizing $2$-sphere $\S$, and then globalize.  A key step in proving the local splitting is to obtain a foliation by constant mean curvature $2$-spheres in a neighborhood of $\S$, as in the proof of  Theorem \ref{rigid}.  A related rigidity result
for least area projective planes has been obtained in \cite{brayetal}.  
Recent counterexamples to Min-Oo's conjecture \cite{brendle} show that ``least area" cannot be replaced by ``stable" in the latter.


\section{Marginally outer trapped surfaces}

Some of the results presented in the previous section can be interpreted as statements concerning
 {\it time-symmetric} initial data sets in general relativity.   In this section we shall consider some  closely related results in general relativity, pertaining to the topology of black holes,
which hold for general initial data sets.  In this more general spacetime setting, minimal surfaces are naturally replaced by the spacetime notion of {\it marginally outer trapped surfaces}.  

The notion of a marginally outer trapped surface  (MOTS) was introduced early on in the development of the theory of black holes.  Under suitable circumstances, the occurrence of a MOTS in a time slice signals the presence of a black hole \cite{HE,CGS}.   For this and other reasons MOTSs have played  a fundamental role in quasi-local descriptions of  black holes; see e.g.,  \cite{AK}.  MOTSs arose in a more purely mathematical context  in the work of Schoen and Yau \cite{SY2} concerning the existence of solutions of Jang's equation, in connection with their proof of positivity of mass.   There have been significant  advances in the analysis of MOTSs in recent years.  We refer the reader to the survey article \cite{AEM} which describes many of these developments, including  the important connections between solutions of Jang's equation and the existence of MOTSs.

Here we are mostly concerned with MOTSs in initial data sets.  
Let $(\bar M^{n+1},\bar  g)$ be a spacetime (time oriented Lorentzian manifold) of dimension $n+1$,
$n \ge 3$.   By an initial data set in $(\bar M^{n+1},\bar  g)$, we mean a
triple $(M^n, g, K)$, where $M$ is a spacelike
hypersurface in $\bar M$, and $g$ and $K$ are the induced metric and second fundamental
form, respectively, of $M$.  To set sign conventions, for vectors $X,Y \in T_pM$, $K$ is defined as, $K(X,Y) = \<\bar \D_X u,Y\>$, where $\bar\D$ is the Levi-Civita connection of $\bar M$ and $u$ is the future directed timelike unit normal vector field to $M$. 

Given an initial data set $(M^n, h, K)$,   let $\S^{n-1}$ be a closed (compact without boundary) two-sided hypersurface in $M^n$.   Then $\S$ admits a smooth unit normal field
$\nu$ in $M$, unique up to sign.  By convention, refer to such a choice as outward pointing. 
Then $l_+ = u+\nu$ (resp. $l_- =  u - \nu$) is a future directed outward (resp., future directed inward) pointing null normal vector field along $\S$, unique up to positive scaling.   

The second fundamental form of $\S$ (viewed as a co-dimension two submanifold of $(\bar M^{n+1},\bar  g)$) can be decomposed into two scalar valued {\it null second forms},  $\chi_+$ and $\chi_-$, associated to $l_+$ and $l_-$, respectively.
For each $p \in \S$, $\chi_{\pm}$ is 
the bilinear form defined by,
\beq
\chi_{\pm} : T_p\S \times T_p\S \to \mathbb R , \qquad \chi_{\pm}(X,Y) = \bar g(\bar\D_Xl_{\pm}, Y) \,.
\eeq
The {\it null expansion scalars} (or {\it null mean curvatures})  $\th_{\pm}$ of $\S$   are obtained by tracing 
$\chi_{\pm}$ with respect to the induced metric $h$ on $\S$,
\begin{align}
\theta_{\pm} = {\rm tr}_{h} \chi_{\pm} = h^{AB}\chi_{\pm AB} = {\rm div}\,_{\S} l_{\pm}  \,.
\end{align}
One verifies that the sign of $\th_{\pm}$ is invariant under positive scaling 
of the null vector field $l_{\pm}$. Physically, $\th_+$ (resp., $\th_-$) measures the divergence
of the  outgoing (resp., ingoing) light rays emanating
from $\S$.  In terms of the initial data $(M^n,h,K)$,
\beq\label{thid}
\th_{\pm} = {\rm tr}_{h} K \pm H \,,
\eeq 
where $H$ is the mean curvature of $\S$ within $M$ (given by the divergence of $\nu$ along~$\S$). 

For round spheres in Euclidean slices of Minkowski space, with
the obvious choice of inside and outside, one has $\th_- < 0$
and $\th_+ >0$. In fact, this is the case in general for large
``radial" spheres in {\it asymptotically flat} spacelike
hypersurfaces.  However, in regions of space-time where the
gravitational field is strong, one may have both $\th_- < 0$
and $\th_+ < 0$, in which case $\S$ is called a {\it trapped
surface}.    Under appropriate energy and causality conditions, 
the occurrence of a trapped surface signals  the onset of gravitational collapse and the existence
of a black hole \cite{HE}.

Focussing attention on the outward null normal only, we say that
$\S$ is an outer trapped surface  if $\th_+ < 0$.  Finally, we define $\S$ to be a marginally
outer trapped surface (MOTS) if $\th_+$ vanishes identically.   

MOTSs arise naturally in a number of situations.  Most basically, compact cross sections of the event horizon in stationary (i.e. steady state) black holes spacetimes  are MOTSs.  This may be understood as follows.  The event horizon is a null hypersurface in spacetime ruled by null geodesics.  By Hawking's area theorem (which requires the null energy condition), these null geodesics can only diverge towards the future.  However in the stationary limit, this divergence vanishes, which implies that cross sections have 
$\th_+ = 0$.   

In dynamical black hole spacetimes MOTSs typically occur inside the event horizon.  There are old heuristic arguments for the existence of MOTSs in this case, based on considering the boundary of the {\it trapped region} inside the event horizon.  These heuristic ideas have recently been
made rigorous first  by Andersson and Metzger \cite{AM2} for three dimensional initial data sets, and then by Eichmair \cite{Eich1, Eich2} for initial data sets up to dimension seven.  These results rely on a basic existence
result for MOTSs under physically natural barrier conditions. The proofs by these authors, based on an approach put forth by  Schoen,  involves inducing blow-up of Jang's equation; 
for further details, see \cite{AEM} and references therein.

Geometrically, MOTSs may be viewed as spacetime analogues of minimal surfaces in Riemannian manifolds.  In fact, in the time-symmetric case (i.e. the case in which $M$ is totally geodesic, $K = 0$), a MOTS $\S$ is just a minimal hypersurface 
in $M$, as can be seen from Equation (\ref{thid}).  Despite the fact that MOTSs are not known to arise
as critical points of some functional,  
they have been shown in recent years to share a number of properties in common with minimal surfaces.  In particular, as first put forward by Andersson, Mars and Simon  \cite{AMS1,AMS2}, MOTSs admit a notion of stability analogous to that of minimal surfaces.  Here, stability is associated with variations of the null expansion under  deformations of a MOTS, as discussed in the following subsection.

\subsection{Stability of MOTSs}

Let $\S$ be a MOTS in the initial data set $(M,g,K)$ with outward unit normal $\nu$.  Consider a normal variation of $\S$ in $M$,  i.e.,  a variation 
$t \to \S_t$ of $\S = \S_0$ with variation vector field 
$V = \frac{\d}{\d t}|_{t=0} = \phi\nu, \quad \phi \in C^{\infty}(\S)$.
Let $\th(t)$ denote the null expansion of $\S_t$
with respect to $l_t = u + \nu_t$, where $u$ is the future
directed timelike unit normal to $M$ and $\nu_t$ is the
outer unit normal  to $\S_t$ in $M$.   A computation shows,
\beq\label{thder} \left . \frac{\d\th}{\d t} \right |_{t=0}   =
\bar L(\f) \;, \eeq where $\bar L : C^{\infty}(\S) \to C^{\infty}(\S)$ is
the operator~\cite{AMS2}, \beq\label{stabop2}
\bar L(\phi)  = -\triangle \phi + 2\<X,\D\phi\>  + \left( \frac12 S_{\S}
- (\mu + J(\nu)) - \frac12 |\chi|^2+{\rm div}\, X - |X|^2
\right)\phi \,. \eeq 
In the above, $\triangle$, $\D$ and ${\rm
div}$ are the Laplacian, gradient and divergence operator,
respectively, on $\S$, $S_{\S}$ is the scalar curvature of $\S$, $X$ is the 
vector field  on $\S$  dual to the one form $K(\nu,\cdot)|_{\S}$, 
$\<\,,\,\> =h$ is the induced metric  on $\S$, and $\mu$ and $J$ are defined in terms
of the Einstein tensor $G = {\rm Ric}_{\bar M} - \frac12 R_{\bar M} \bar g$\,: $\mu = G(u,u)$,
$J = G(u,\cdot)$.  When the Einstein equations are assume to hold (see Section 2.2), $\mu$ and $J$ represent the  energy density and linear momentum density along $M$.  As a consequence of the Gauss-Codazzi equations, the quantities $\mu$ and $J$ can be expressed solely in terms of initial
data,
\beq\label{emid}
\mu = \frac12\left(S_{\S} + ({\rm tr}\,K)^2 - |K|^2\right) \quad \text{and} \quad
J = \div K- d({\rm tr}\, K)  \,.
\eeq

In the time-symmetric case,  $\th$ in \eqref{thder} becomes the
mean curvature $H$, the vector field $X$ vanishes and $\bar L$
reduces to the stability operator \eqref{stabop} of minimal surface theory. In
analogy with the minimal surface case, we refer to $\bar L$ in
\eqref{stabop2} as the stability operator associated with
variations in the null expansion $\th$.  As observed in \cite{AMS2}, although in general
$\bar L$ is not self-adjoint, as a consequence of the Krein-Rutman theorem, its principal eigenvalue (eigenvalue with smallest real part) $\l_1(\bar L)$ is real.   Moreover there
exists an associated eigenfunction $\phi$ which is strictly positive.

Continuing the analogy with the minimal surface case, we
say that a MOTS is stable provided $\l_1(\bar L) \ge 0$.   Stable MOTSs arise in various
situations.   For example, outermost MOTSs, as we define in Section 2.2, are stable.
This includes, in particular, compact cross sections of the event horizon in stationary black
hole spacetimes obeying the null energy condition.  More generally, the results of Andersson and
Metzger, and of Eichmair alluded to above establish natural criteria for the existence of 
outermost MOTSs.  And of course, in the time-symmetric case, stable minimal surfaces are stable
MOTSs. 

We now present a key fact which allows, for example, some of the results on minimal surfaces in Section 1 to be extended to MOTSs in general initial data sets.   Consider the
``symmetrized" operator
$L_0: C^{\infty}(\S) \to C^{\infty}(\S)$,
\beq\label{symop}
L_0(\phi)  = -\triangle \phi  + \left( \frac12 S_{\S} - (\mu + J(\nu)) - \frac12 |\chi|^2\right)\phi \,.
\eeq
obtained by formally  setting $X= 0$ in \eqref{stabop2}.   The key argument in~\cite{GS}
used to obtain a generalization of Hawking's black hole topology theorem \cite{HE}
shows the following.

\begin{prop}\label{eigen}
$\l_1(L_0) \ge \l_1(\bar L)$.
\end{prop}

\proof  The proof has its roots in arguments from \cite{SY2}.  Let $\phi$ be a positive eigenfunction associated to $\l_1(\bar L)$.  Completing the square on the right hand side of 
\eqref{stabop2},  and using $\bar L(\phi) = \l_1(\bar L)\phi$ gives,
\beq
 -\triangle \phi +\left(Q+{\rm div}\, X  \right)\phi   +
\phi|\D \ln\phi|^2  - \phi|X - \D\ln\phi|^2  = \l_1(\bar L) \phi
\eeq
where  $Q:=  \frac12 S - (\mu + J(\nu)) - \frac12 |\chi|^2$.
Letting $u = \ln \phi$, we obtain,
\beq\label{uineq}
-\triangle u +Q + {\rm div}\,X  - |X - \D u|^2  =  \l_1(\bar L)  \,.
\eeq
Absorbing the Laplacian term
$\triangle u = {\rm div}\,(\D u)$  into the divergence term 
gives,
 \beq\label{stabid}
Q + {\rm div}\,(X- \D u)  - |X - \D u|^2 =  \l_1(\bar L).
\eeq
Setting $Y = X - \D u$, we arrive at,
\beq
- Q + |Y|^2  +  \l_1(L) =  {\rm div}\, Y  \,.
\eeq

Given any $\psi \in C^{\infty}(\S)$, we multiply  through by $\psi^2$ and derive,
\begin{align*}
-\psi^2 Q +\psi^2 |Y|^2 + \psi^2\l_1(L)  &= \psi^2 {\rm div}\, Y \nonumber\\
& =  {\rm div}\,(\psi^2Y) - 2\psi \< \D\psi,Y \>  \nonumber  \\
& \le   {\rm div}\,(\psi^2Y) + 2|\psi| |\D \psi| |Y| \nonumber\\
& \le  {\rm div}\,(\psi^2Y) + |\D \psi|^2 + \psi^2|Y|^2  \, .
\end{align*}

Integrating the above inequality yields,
\beq
 \l_1(\bar L) \le \frac{\int_{\S} |\D \psi|^2 + Q \psi^2}{\int_{\S} \psi^2}  \quad \mbox{for all } \psi \in C^{\infty}(\S), \, \psi \not\equiv 0  \,.
\eeq
Proposition \ref{eigen} now follows from the well-known Rayleigh formula for the principal eigenvalue
applied to the operator \eqref{symop}.~\qed

\medskip
When $\S$ is stable, the left hand side of (\ref{stabid}) is nonnegative.  The resulting inequality is intimately related to the  {\it stability inequality}  (2.29) for solutions to Jang's equation obtained by Schoen and Yau in \cite{SY2}.  
For further discussion of the connection between Jang's equation and MOTSs, see e.g. \cite{AEM}.

\subsection{A  generalization of Hawking's black hole topology theorem}

A basic step in the proof of the classical
black hole uniqueness theorems is Hawking's
theorem on the topology of black holes~\cite{HE} which asserts
that compact cross-sections of the event horizon in
$3+1$-dimensional asymptotically flat
stationary black hole space-times obeying the dominant energy
condition are topologically 2-spheres.   The remarkable discovery of
Emparan and Reall~\cite{ER} of a $4+1$ dimensional asymptotically flat
stationary vacuum black hole space-time with horizon topology $S^1
\times S^2$, the so-called ``black ring", showed that horizon
topology need not be spherical in higher dimensions.  This example naturally led to the
question of what are the allowable horizon topologies in higher
dimensional black hole space-times. This question was addressed
in the papers of~\cite{GS, G}, resulting
in a natural generalization of Hawking's black hole topology theorem to
higher dimensions, which we now discuss.

Let $(\bar M^{n+1}, \bar g)$, $n \ge 3$, be a spacetime satisfying the Einstein equations
(with vanishing cosmological term),
\beq
{\rm Ric}_{\bar M} - \frac12 R_{\bar M} \bar g = \cal T \,,
\eeq
where ${\rm Ric}_{\bar M}$ and $ R_{\bar M}$ are the Ricci tensor and scalar curvature of   
$(\bar M^{n+1}, \bar g)$, respectively, and  $\cal T$ is the energy momentum tensor.  
$(\bar M^{n+1}, \bar g)$ is said to satisfy the dominant energy condition (DEC) provided
$\calT(X,Y) = T_{ij}X^iY^j \ge 0$ for all future directed causal vectors $X,Y$.  This translates to
the spacetime curvature condition,  $G(X,Y) = G_{ij}X^iY^j \ge 0$, where $G = {\rm Ric}_{\bar M} - \frac12 R_{\bar M} \bar g$ is the Einstein tensor.   One verifies that if $(\bar M^{n+1}, \bar g)$
satisfies the DEC then for every initial data set $(M^n,g, K)$ in  $(\bar M^{n+1}, \bar g)$ one has
$\mu \ge |J|$ along $M$, where $\mu = G(u,u)$ and $J = G(u, \cdot)$.  Recall from equation
\eqref{emid}, that the energy density $\mu$ and linear momentum density $J$ along $M$ can be
expressed solely in terms of initial data.  Note that in the time-symmetric case, $\mu = \frac12 S_{\S}$ and $J = 0$. 

The following result of Schoen and the author \cite{GS}, which generalizes Hawking's black hole topology theorem, is a spacetime analogue of Proposition \ref{pos} in Section~1.
\begin{thm}\label{pos2} 
Let  $(M^n,h,K)$, $n \ge~3$, be an initial data set such 
$\mu \ge |J|$ along $M$.
If $\S$ is a stable MOTS in $M$ then $\S$ is of positive Yamabe type, unless
$\S$ is Ricci flat (flat if $n=3,4$), $\chi = 0$ and $\mu + J(\nu) = 0$ along $\S$.
\end{thm} 

\proof[Comment on the proof.]  Since $\S$ is stable, $\l_1(\bar L) \ge 0$, where $\bar L$ is the 
MOTS stability operator \eqref{stabop2}.  By Proposition \ref{eigen}, we have $\l_1(L_0) \ge 0$.
Noting the similarity between $L_0$ and the minimal surface stability operator $L$ in \eqref{stabop}, we see that to complete the proof of Theorem \ref{pos2}, one can argue just as in the proof of Proposition \ref{pos}, by using a positive eigenfunction $\phi$ associated to $\l_1(L_0)$ to construct the relevant conformally related metric.
The only other difference is that the nonnegative quantity $\frac12 S + \frac12 |B|^2$ appearing 
in \eqref{eq:conf} is now replaced by the nonnegative quantity $\mu + J(\nu) + \frac12 |\chi|^2$.
\qed

\medskip
Thus, apart from certain exceptional circumstances, a stable MOTS $\S$ in an initial data set satisfying the DEC, must be of positive Yamabe type.  But $\S$ being positive Yamabe implies many well-known restrictions on its topology.  We single out a couple of cases here; for simplicity assume $\S$
is orientable.  In the standard case: dim $\S = 2$ (dim $M = 3+1$), $\S$ admits a metric of positive Gaussian curvature and so is diffeomorphic to a $2$-sphere by Gauss-Bonnet.  Hence we recover
Hawking's theorem.  In the case  dim $\S = 3$ (dim $M = 4+1$), it is known by results of Schoen-Yau \cite{SYsc2} and Gromov-Lawson \cite{GL} that $\S$ must be diffeomorphic to (i) a spherical 
space (a  quotient of the $3$-sphere - by the positive resolution of the Poincar\'e conjecture!), or
(ii) $S^2 \times S^1$ or (iii) a connected sum of the previous two types.  Thus, the basic horizon topologies in the case dim $\S = 3$ are $S^3$ and $S^2 \times S^1$, the latter being realized
by the black ring of Emparan and Reall.  For further restrictions on horizon topology within the class of stationary black hole spacetimes in $4+1$ dimensions, see \cite{Hollands}.

A drawback of Theorem \ref{pos2} is that, when the DEC along $(M,g,K)$ does not hold strictly, it allows
certain possibilities that one would like to rule out.  For example, it does not rule out the possibility
of a vacuum black hole spacetime with toroidal topology.  In fact it is easy to construct initial data
sets obeying the DEC with stable toroidal MOTSs.   However these exceptional cases 
can be eliminated by strengthening the stability assumption.  

As noted in Section 2, stable MOTSs often arise as {\it outermost MOTSs}. 
Let $\S$ be a  MOTS  in an initial data set $(M,g,K)$ that separates $M$ into an ``inside" and
an ``outside". We say that $\S$ is  outermost if
there are no outer trapped ($\th_+  < 0$)  or marginally outer trapped ($\th_+= 0$) 
surfaces outside of and homologous to $\S$. If $\S$ is outermost it is necessarily stable,
for otherwise, using a positive eigenfunction $\phi$ associated to $\l_1(\bar L) < 0 $ to define a 
normal variation $t \to \S_t$ of $\S = \S_0$, one sees that for small $t >0$, $\S_t$ is outer trapped.
(The term ``outermost" here is partially justified by the fact that, under a natural exterior barrier
condition, and up to dimension seven, if there exists an outer trapped surface $\S'$ outside of 
and homologous to $\S$, then there will exist a MOTS outside of and homologous to $\S'$; see
\cite{AEM} and references cited therein.)

The following result, roughly similar to Corollary \ref{outermost}, and proved in \cite{G}, rules out
the exceptional cases. 

\begin{thm}
Let $(M^n,g,K)$, $n \ge 3$, be an initial data set in a spacetime obeying the DEC. If $\S^{n-1}$
is an outermost  MOTS in $(M^n,g,K)$ then it is of positive Yamabe type. 
\end{thm}

This theorem is an immediate consequence of the following rigidity result (compare Theorem 
\ref{rigid}).

\begin{thm}\label{rigid2}
Let $(M^n,h,K)$, $n \ge 3$, be an initial data set such that $\mu \ge |J|$ along $M$, and such that
$\tau = {\rm tr}_g K \le 0$.  Suppose  $\S$ is a separating MOTS in $M$  such that there are
no outer trapped surfaces ($\th_+ < 0$) outside of, and homologous, to $\S$.  If $\S$ 
is  not of positive Yamabe type,  then there exists
an outer half-neighborhood $U \approx [0,\e) \times \S$ of $\S$ in $M$ such that each slice
$\S_t = \{t\} \times \S$, $t \in [0,\e)$ is a MOTS. In fact each $\S_t$ has $\chi_t = 0$, and is Ricci flat.   Moreover, if the DEC 
holds in a spacetime neighborhood of $\S$ then the mean curvature condition, $\tau \le 0$ can be dropped.
\end{thm}

\proof[Remarks on the proof.]  The proof of  Theorem \ref{rigid2}, while somewhat more involved, is similar to the proof given here of Theorem \ref{rigid}.   First one uses stability  to obtain
an outer foliation $t \to \S_t$, $0 \le t \le \e$, of surfaces $\S_t$ of {\it constant} outer null expansion,
 $\th(t) = c_t$.  Then one shows that the constants $c_t = 0$.  The proof of this makes use of the formula for the $t$-derivative, $\frac{d\th}{d t}$, not just at $t=0$ where
$\th =0$, but all along the foliation $t \to \S_t$, where, a priori, $\th(t)$ need not be zero.  Thus,  additional terms appear in the expression for $\frac{d \th}{d t}$ beyond those appearing in (\ref{thder})-({\ref{stabop2}), including a term involving the mean curvature $\tau$ of $M$, which need to be accounted for.  By making a small deformation of $M$ near $\S$ in spacetime, one can eliminate
the mean curvature condition on $M$, provided the DEC holds in a neighborhood of $\S$.  We are working on an approach to eliminate this deformation step, so as to obtain a pure initial data version
of Theorem \ref{rigid2}, without the mean curvature condition.

\section{Positivity of mass for asymptotically hyperbolic manifolds}

Without a doubt one of the most beautiful and powerful applications of minimal surface theory to manifolds
of nonnegative scalar curvature is Schoen and Yau's proof of the positive mass theorem \cite{SY1} for maximal
three-dimensional initial data sets.
In \cite{S} Schoen  showed by an inductive argument how to extend their positive mass result up to dimension seven.

\begin{thm}[Riemannian PMT]
Suppose $(M^n,g)$, $3\le n\le 7$,  is an asymptotically flat 
manifold with
nonnegative scalar curvature, $S\ge 0$.  Then $M$ has ADM mass $\ge 0$, and $=0$ iff
M is isometric to Euclidean space.
\end{thm}


Developments in physics over the past decade, in particular the emergence of the AdS/CFT correspondence,
has heightened interest in the geometrical and physical properties of asymptotically hyperbolic (AH) manifolds.
Such manifolds arise naturally as spacelike hypersurfaces in asymptotically anti-de Sitter spacetimes. 
In \cite{Gibbons}, Gibbons et al. adapted Witten's spinorial argument to prove positivity of mass in the $3+1$ dimensional asymptotically AdS setting.
More recently, Wang \cite{Wang}, and, under  weaker asymptotic conditions,
Chru\'sciel and Herzlich~\cite{ChH} provided precise definitions of the mass in the asymptotically hyperbolic setting and gave spinor based proofs of positivity of mass in dimensions $\ge 3$.  These latter positive mass results may be
paraphrased  as follows:

\begin{thm}\label{pmass}
Suppose $(M^{n},g)$, $n \ge 3$, is an asymptotically hyperbolic spin manifold with
scalar curvature $S\ge -n(n-1)$.  Then $M$ has mass $\ge 0$, and $=0$ iff M
is isometric to standard hyperbolic space $\mathbb H^{n}$.
\end{thm}

Physically, $M$ corresponds to a maximal (mean curvature zero)  spacelike
hypersurface in spacetime satisfying the Einstein equations with cosmological
constant $\Lambda = -n(n-1)/2$.  For then the Gauss equation and DEC
imply $S \ge -n(n-1)$.  We note that, in the present context, the mass  in Theorem \ref{pmass} refers, more precisely,
to  the energy component of a suitably defined energy-momentum vector \cite{ChH, Wang}.

The aim of the paper \cite{ACG} was to obtain a positive mass theorem
for asymptotically hyperbolic (AH) manifolds which does not require the assumption of spin, with
the idea of replacing the spinor based arguments by the ``minimal surface"
methodology of Schoen and Yau, adapted to a negative lower bound on the scalar curvature.
The specific approach taken in \cite{ACG} was
inspired by Lohkamp's \cite{Lohk} variation of the Schoen-Yau proof.  Lohkamp 
circumvented many of the technical issues in \cite{SY1} associated with the construction 
and analysis of a complete noncompact  minimal surface, by {\it compactifying} the
problem.  This allows one to appeal to results on compact manifolds of nonnegative scalar curvature in
\cite{SYsc2} which do not require the assumption of spin.   Recall, under the assumptions
of Theorem \ref{pmass}, Lohkamp shows that if the mass is negative then the metric $g$ can
be deformed to a metric $\tilde g$ having nonnegative scalar curvature, $\tilde S \ge 0$ (with 
$\tilde S > 0$ at some points) such that $(M,\tilde g)$ is isometric to Euclidean space $(\Bbb R^n,g_0)$ 
outside a compact set.  But then  this contradicts the following result, first proved by Schoen
and Yau in dimension three \cite{SYNAS, SYsc}.

\begin{thm}[Precursor to PMT]\label{precursor}
Suppose that ($M^n,g)$, $3 \le n \le 7$,  has nonnegative scalar curvature, $S\ge 0$, 
and is isometric to $(\Bbb R^n,g_0)$ outside a compact set $K$.  Then $(M^n,g)$ is {\it globally} isometric to $(\Bbb R^n,g_0)$.
\end{thm}

Recall, to prove Theorem \ref{precursor}, one encloses $K$ in a large $n$-dimensional cube, with sides
in the Euclidean region.  Then, by identifying sides pairwise one obtains a compact manifold  of the form
$T^n \# N$ with nonnegative scalar curvature.  By Corollary 2 in \cite{SYsc2}, such a manifold must be flat.  (If $N$ is spin, this follows from results of Gromov and Lawson for all  $n \ge 3$; see \cite{GL} and references cited therein.)
Thus, $(M^n,g)$ is globally flat, and Theorem \ref{precursor}  follows.  

The idea of  the paper \cite{ACG} was to carry out a similar plan  to prove a version of Theorem \ref{pmass}
without spin assumption.  This plan consists of two parts: a deformation part and a rigidity part.  The defomation part  is to show in the AH setting, that if the mass is negative then one can
deform the metric near infinity to make it exactly hyperbolic outside an arbitrarily large compact
set, while retaining the scalar curvature inequality, $S \ge -n(n-1)$.  We will say a bit more about  this part at the end.  The rigidity part is to obtain a hyperbolic analogue of Theorem~\ref{precursor}.   In fact, in \cite{ACG} the following rigidity result was established.

\begin{thm}\label{Hrigid}
Suppose $(M^{n},g)$, $3 \le n \le 7$, has scalar curvature $S$ 
satisfying, $S\ge -n(n-1)$,
and is isometric to $\mathbb H^{n}$  outside a compact set.
Then $(M^{n},g)$ is globally isometric to $\mathbb H^{n}$.
\end{thm}

In the case that $(M^{n},g)$ is a spin manifold, this theorem 
follows from a result of Min-oo \cite{Minoo} (see also \cite{AD,Delay}), as well as from  the
rigidity part of Theorem \ref{pmass}. The main point of
Theorem \ref{Hrigid} is that it does not require   a spin assumption.
We wish to take some time here to discuss the proof, as it relates closely to
themes already developed in the first two sections. The proof is based on a study of minimizers of 
the `BPS brane action', as it is referred to in \cite{WY}.

\subsection{The brane action}

Let $(M^{n}, g)$ be an oriented Riemannian manifold
with volume form $\Omega$.  Assume there is a globally defined
form $\Lambda$ such that $\Omega = d \Lambda$.
Let $\S$ be a compact orientable hypersurface in $M$.   Then $\S$ is
$2$-sided in $M$, and hence admits a smooth unit normal field $\nu$, which we
refer to as outward pointing.   Let $\S$ have the orientation induced by $\nu$ (i.e., determined by
the induced area form $\omega = i_{\nu}\Omega$).  
Then, for any such
$\S$, we define the brane action $\calB$ by,
\beq
\calB(\S) = \calA(\S) -(n-1) \calV(\S) \, ,
\eeq
where $\calA(\S) =$ the area of $\S$, and $\calV(\S) =  \int_{\S} \Lambda$.
If $\S$ bounds to the inside then,  by Stokes theorem, $\calV(\S) =$ the volume of the region enclosed by $\S$.   Although $\Lambda$ is not uniquely determined,
Stokes theorem shows that, within a given homology class, $\calB$ is uniquely
determined up to an additive constant.

We now consider the formulas for the first and second variation of the
brane action.  Let $t \to \S_t$, $-\e < t < \e$, be a normal variation of $\S= \S_0$, with variation vector
field $V = \left . \frac{\d}{\d t}\right |_{t=0} = \phi \nu$,  $\phi \in C^{\infty}(\S)$.   Abusing notation slightly, set
$\calB(t) = \calB(\S_t)$.  Then for first variation we have,
\beq
\calB'(0) = \int_{\S}  (H -(n-1))\f \, dA
\eeq
Thus $\S$ is a stationary point for the brane action if and only if it
has constant mean curvature $H = n-1$.

Assuming $\S$ has mean curvature
$H =n-1$, the second variation formula is given by
\beq
\calB''(0) = \int_{\S} \phi L(\phi)\, dA  \,,
\eeq
where $L$ is just as in \eqref{stabop0}.
Using \eqref{gauss} and the fact that $H =n-1$,  $L$ can be expressed as,
\beq\label{stabop3}
L(\phi) = - \L \phi + \frac12(S_{\S} - S_n - |B_0|^2)\,\phi   \,,
\eeq
where $S_n = S + n(n-1)$ and $B_0$ is the trace free part of $B$, $B_0 = B -h$, where $h$
is the induced metric on $\S$.
We note that, in our applications, $S_n$ will be nonnegative.

A stationary point $\S$ for the brane action is said to be $\calB$-stable provided
for all normal variations $t \to \S_t$ of $\S$, $\calB''(0) \ge 0$. As in Section 1,
$\S$ is  $\calB$-stable if and only if $\l_1(L) \ge 0$, where $L$ is the brane stability
operator \eqref{stabop3}.

In close analogy with Proposition \ref{pos}, one obtains, by essentially the same argument, the following.
\begin{prop}\label{pos3}
Let $(M^{n},g)$, $n \ge 3$, be an oriented Riemannian manifold having scalar curvature $S \ge -n(n-1)$. 
If $\S$ is a compact orientable
$\calB$-stable hypersurface in $M$ which does not admit a metric of positive
scalar curvature then (i)~$\S$ is umbilic; more precisely, $B =h$, where $h$ is the induced metric on 
$\S$, and (ii)~$\S$ is Ricci flat and $S = -n(n-1)$ along $\S$.
\end{prop}

This infinitesimal rigidity result is used in the proof of the following warped product splitting result.  

\begin{thm}\label{warp}
Let ($M^{n},g)$ be an oriented Riemannian manifold
with scalar curvature $S \ge -n(n-1)$.  Let $\S$ be a compact orientable
hypersurface in $M$ which does not admit a metric of positive scalar curvature.
If $\S$ locally minimizes the brane action $\calB$ then there is a neighborhood
$U$ of $\S$ such that $(U,g|_U)$ is isometric to the warped product
$((-\e,\e)\times \S, dt^2 +e^{2t} h)$, where $h$, the induced metric on $\S$, is
Ricci flat.
\end{thm}

As noted in Section 1, the proof of  Theorem \ref{warp} given in \cite{ACG} is very similar to the proof given here of Theorem \ref{rigid}.  A related result in three dimensions 
has been obtained in \cite[Theorem 3.2]{Yau}. 

\subsection{Proof of Theorem \ref{Hrigid}}

We briefly outline the proof of Theorem \ref{Hrigid}.  While, unlike the asymptotically flat case,
there  does not appear to be a way to fully compactify, our approach taken in \cite{ACG} is 
to {\it partially compactify} and then minimize the brane action in a certain homology class.  

We are assuming that outside a compact set $K$, $(M^n,g)$ is isometric to hyperbolic space $\bbH^n$.
For our purposes, it is convenient to work with an explicit representation of  $\bbH^{n}$.  We start with
the half-space model $(H^{n},g_H)$, where, $H^{n} = \{(y, x^1,\cdots, x^{n-1}) : y > 0\}$, and
\beq
g_H = \frac1{y^2} \left (dy^2 +  (dx^1)^2 + \cdots + (dx^{n-1})^2 \right)  \,,
\eeq
and make the change of variable $y= e^{-t}$, to obtain
$\bbH^{n}= (\bbR^{n}, g_1)$,
where,
\beq
g_1 = dt^2 + e^{2t} \left( (dx^1)^2 + \cdots + (dx^{n-1})^2 \right)  \,.
\eeq

To partially compactify, we periodically identify the $x^i$ coordinates.  Choosing $a>0$ sufficiently large, we
can enclose the compact set $K$ in an infinitely long rectangular box, with sides
determined by the ``planes",
$x^i = \pm a$.  We  then identify the sides of the box pairwise  to obtain an identification space
$(\hat M, \hat g)$ such that outside the compact set $K$, and up to isometry we have,
\beq\label{cusp}
\hat M = \bbR \times T^{n-1} \,,  \quad \hat g = dt^2 + e^{2t} h \,,
\eeq
where $h$ is a flat metric on the torus $T^{n-1}$.  Thus, $(\hat M, \hat g)$ is just
a standard hyperbolic cusp ``perturbed" on the compact set $K$, 
with scalar curvature satisfying $\hat S \ge -n(n-1)$
everywhere.

\medskip
\noindent
{\bf Claim.} $(\hat M, \hat g)$ has constant curvature $-1$ everywhere.


\medskip

The claim implies that the original manifold $(M,g)$, from which $(\hat M, \hat g)$ was constructed,  has constant curvature $-1$ everywhere. This implies that $(M,g)$ is covered by $\bbH^n$, with the covering map a local isometry.  But since $(M,g)$ is simply connected
at infinity, the covering map must be an isometry, and Theorem~\ref{Hrigid} follows.

We give a sketch of the proof of the claim.
Choose $b > 0$ large enough so that $K$ lies between the
toroidal slices $t = \pm b$, and fix
a $t$-slice $\S_{t_0} = \{t_0\} \times T^n$, $t_0 > b$, with outward pointing normal $\nu = \frac{\d}{\d t}$.  We now minimize the brane action in the homology class determined by $\S_{t_0}$.
Let $S_i$ be a minimizing sequence for the brane action in this homology class.  The $t$-slices $\S_t$, $|t|  \ge b$, have mean curvature $H =n-1$,
and hence act as  barriers for the minimization procedure. Consequently, the $S_i's$ can be chosen to lie in the
compact region between the slices $t= \pm b$. Thus, by the compactness results of GMT, a subsequence
converges to a smooth (provided $n \le 7$)  compact hypersurface $S$ homologous to $\S_{t_0}$, which minimizes the brane action in its
homology class. Without loss of generality we may assume $S$ is connected; otherwise
restrict attention to one of its components.  The ``almost product structure" of $\hat M$ (cf., Equation \eqref{cusp}) implies the existence of a nonzero degree map from $S$ to the toroidal slice $\S_{t_0}$.  
By Corollary 2 in \cite{SYsc}, $S$ does not admit a metric of positive scalar curvature (and any metric of nonnegative scalar curvature on $S$ must be flat).  
Then, by Theorem \ref{warp}, a neighborhood $U$ of $S$ splits as a warped
product,
\beq
U = (-u_0, u_0) \times S \qquad \hat g|_U = du^2 + e^{2u} h  \, ,
\eeq
where the induced metric $h$ on $S$ is flat.  But since $S$  globally
minimizes the brane action in its homology class,  this local
warped product structure can be extended to arbitrarily large $u$-intervals.  Hence
$K$ will eventually be contained in this constructed warp product region. The claim now follows.
(In fact, with a bit more effort, these arguments show that Equations \ref{cusp} hold globally.)  \qed

\medskip
We conclude with a few words about the deformation step.  In \cite{ACG}, we adopted the definition
of the mass of an asymptotically hyperbolic manifold used by Wang~\cite{Wang}. By a change of coordinates
this implies that there is a relatively compact set $C$
such that   $M\setminus C = S^n \times [R,\infty)$, $R > 0$, 
and on $M \setminus C$, $g$ has  the form,
\beq\label{metric}
g = \frac1{1+r^2} \, dr^2 + r^2h_r \,, 
\eeq
where $h_r$ is an $r$-dependent family of metrics on $S^n$ of the form,
\beq
h_r = h_0+ \frac1{r^{n}}\, k + \g_r \,,
\eeq
where $h_0$ is the standard metric on $S^n$, $k$ is a symmetric $2$-tensor on $S^n$ and $\g_r$ is an $r$-dependent family of metrics on $S^n$ whose components  satisfy specific higher order decay rates.
We refer to $k$ as the mass aspect tensor; it is the leading order measure of the deviation of the metric $g$
from the hyperbolic metric.   Its trace with respect to $h_0$, 
${\rm tr}_{h_0}\,k$,
is called the mass aspect function.
Up to a normalizing constant, the integral of the mass aspect function 
over the sphere defines the mass of $(M,g)$, mass $= \int_{S^n} {\rm tr}_{h_0}\, k$.

What is shown in \cite{ACG}, is that if the mass aspect function is pointwise negative, then sufficiently  far out
on the end  we can deform the metric to the hyperbolic metric, while preserving the scalar curvature inequality $S \ge -n(n-1)$; see Theorem 3.2 in \cite{ACG} for a more precise statement.  But this deformation result is seen to be incompatible with Theorem \ref{Hrigid}.  Moreover, using
Theorem \ref{Hrigid} it is shown  that if the mass aspect function vanishes identically then $(M^n,g)$ is isometric to hyperbolic space.  Putting all of these results together yields the following positive mass statement obtained in~\cite{ACG}.

\begin{thm} \label{pmass2} 
Let $(M^{n},g)$, $3 \leq n \leq 7$, be an asymptotically hyperbolic
manifold with scalar curvature $S \ge -n(n+1)$. Assume that the mass
aspect function does not change sign, i.e. that it is everywhere either negative,
zero, or positive. Then, either the mass of $(M,g)$ is positive, or $(M,g)$
is isometric to hyperbolic space. 
\end{thm}

Naturally, it would be desirable to find a way to remove the sign condition
on the mass aspect.
Within the context of the approach taken in \cite{ACG}, one possible way
to accomplish this would be to suitably extend the results of Corvino-Schoen \cite{CS}
and Chru\'sciel-Delay \cite{CD} on initial data
deformations to the asymptotically hyperbolic setting.   Some partial progress in this direction
has been made in \cite{CD2}.

\providecommand{\bysame}{\leavevmode\hbox to3em{\hrulefill}\thinspace}


\begin{thebibliography}{99}

\bibitem{ACG}
L.~Andersson, M.~Cai, and G.~J. Galloway, \emph{Rigidity and positivity of mass
  for asymptotically hyperbolic manifolds}, Ann. Henri Poincar\'e \textbf{9}
  (2008), no.~1, 1--33.

\bibitem{AD}
L.~Andersson and M.~Dahl, \emph{Scalar curvature rigidity for asymptotically
  locally hyperbolic manifolds}, Ann. Global Anal. Geom. \textbf{16} (1998),
  no.~1, 1--27.

\bibitem{AEM}
L.~Andersson, M.~Eichmair, and J.~Metzger, \emph{Jang's equation and its
  applications to marginally trapped surfaces},  arXiv:1006.4601,
  to appear in the Proceedings of the Conference on Complex Analysis and Dynamical Systems
  IV, Nahariya, 2009.

\bibitem{AMS1}
L.~Andersson, M.~Mars, and W.~Simon, \emph{Local existence of dynamical and
  trapping horizons}, Phys. Rev. Lett. \textbf{95} (2005), 111102.

\bibitem{AMS2}
\bysame, \emph{Stability of marginally outer trapped surfaces and existence of
  marginally outer trapped tubes}, Adv. Theor. Math. Phys. \textbf{12} (2008),
  no.~4, 853--888.

\bibitem{AM2}
L.~Andersson and J.~Metzger, \emph{The area of horizons and the trapped
  region}, Comm. Math. Phys. \textbf{290} (2009), no.~3, 941--972.

\bibitem{AK}
Abhay Ashtekar and Badri Krishnan, \emph{Isolated and dynamical horizons and
  their applications}, Living Reviews in Relativity \textbf{7} (2004), no.~10.

\bibitem{brayetal}
H.~Bray, S.~Brendle, M.~Eichmair, and A.~Neves, \emph{Area-minimizing
  projective planes in three-manifolds},  arXiv:0909.1665.

\bibitem{bray}
H.~Bray, S.~Brendle, and A.~Neves, \emph{Rigidity of area-minimizing
  two-spheres in three-manifolds},  arXiv:1002.2814.

\bibitem{brendle}
S.~Brendle, F.~C. Marques, and A.~Neves, \emph{Deformations of the hemisphere
  that increase scalar curvature},  arXiv:1004.3088.

\bibitem{Cai}
M.~Cai, \emph{Volume minimizing hypersurfaces in manifolds of nonnegative
  scalar curvature}, Minimal surfaces, geometric analysis and symplectic
  geometry ({B}altimore, {MD}, 1999), Adv. Stud. Pure Math., vol.~34, Math.
  Soc. Japan, Tokyo, 2002, pp.~1--7.

\bibitem{CG}
M.~Cai and G.~J. Galloway, \emph{Rigidity of area minimizing tori in
  3-manifolds of nonnegative scalar curvature}, Comm. Anal. Geom. \textbf{8}
  (2000), no.~3, 565--573.

\bibitem{CD}
P.~T. Chru{\'s}ciel and E.~Delay, \emph{On mapping properties of the general
  relativistic constraints operator in weighted function spaces, with
  applications}, M\'em. Soc. Math. Fr. (N.S.) (2003), no.~94, vi+103.

\bibitem{CD2}
\bysame, \emph{Gluing constructions for asymptotically hyperbolic manifolds
  with constant scalar curvature}, Comm. Anal. Geom. \textbf{17} (2009), no.~2,
  343--381.

\bibitem{CGS}
P.~T. Chru{\'s}ciel, G.~J. Galloway, and D.~Solis, \emph{Topological censorship
  for {K}aluza-{K}lein space-times}, Ann. Henri Poincar\'e \textbf{10} (2009),
  no.~5, 893--912.

\bibitem{ChH}
P.~T. Chru{\'s}ciel and M.~Herzlich, \emph{The mass of asymptotically
  hyperbolic {R}iemannian manifolds}, Pacific J. Math. \textbf{212} (2003),
  no.~2, 231--264.

\bibitem{CS}
J.~Corvino and R.~Schoen, \emph{On the asymptotics for the vacuum {E}instein
  constraint equations}, J. Differential Geom. \textbf{73} (2006), no.~2,
  185--217.

\bibitem{Delay}
E.~Delay, \emph{Analyse pr\'ecis\'ee d'\'equations semi-lin\'eaires elliptiques
  sur l'espace hyperbolique et application \`a la courbure scalaire conforme},
  Bull. Soc. Math. France \textbf{125} (1997), no.~3, 345--381.

\bibitem{Eich1}
M.~Eichmair, \emph{The {P}lateau problem for marginally outer trapped
  surfaces}, J. Differential Geom. \textbf{83} (2009), no.~3, 551--583.

\bibitem{Eich2}
\bysame, \emph{Existence, regularity, and properties of generalized apparent
  horizons}, Comm. Math. Phys. \textbf{294} (2010), no.~3, 745--760.

\bibitem{ER}
R.~Emparan and H.~S. Reall, \emph{A rotating black ring solution in five
  dimensions}, Phys. Rev. Lett. \textbf{88} (2002), no.~10, 101101, 4.

\bibitem{FCS}
D.~Fischer-Colbrie and R.~Schoen, \emph{The structure of complete stable
  minimal surfaces in {$3$}-manifolds of nonnegative scalar curvature}, Comm.
  Pure Appl. Math. \textbf{33} (1980), no.~2, 199--211.

\bibitem{G}
G.~J. Galloway, \emph{Rigidity of marginally trapped surfaces and the topology
  of black holes}, Comm. Anal. Geom. \textbf{16} (2008), no.~1, 217--229.

\bibitem{GS}
G.~J. Galloway and R.~Schoen, \emph{A generalization of {H}awking's black hole
  topology theorem to higher dimensions}, Comm. Math. Phys. \textbf{266}
  (2006), no.~2, 571--576.

\bibitem{Gibbons}
G.~W. Gibbons, S.~W. Hawking, G.T. Horowitz, and M.~J. Perry, \emph{Positive
  mass theorems for black holes}, Comm. Math. Phys. \textbf{88} (1983), no.~3,
  295--308.

\bibitem{GL}
M.~Gromov and H.~B. Lawson, Jr., \emph{Positive scalar curvature and the
  {D}irac operator on complete {R}iemannian manifolds}, Inst. Hautes \'Etudes
  Sci. Publ. Math. (1983), no.~58, 83--196 (1984).

\bibitem{HE}
S.~W. Hawking and G.~F.~R. Ellis, \emph{The large scale structure of
  space-time}, Cambridge University Press, London, 1973, Cambridge Monographs
  on Mathematical Physics, No. 1.

\bibitem{Hitchin}
N.~Hitchin, \emph{Harmonic spinors}, Advances in Math. \textbf{14} (1974),
  1--55.

\bibitem{Hollands}
S.~Hollands, J.~Holland, and A.~Ishibashi, \emph{Further restrictions on the
  topology of stationary black holes in five dimensions}, arXiv:1002.0490.

\bibitem{HI}
G.~Huisken and T.~Ilmanen, \emph{The inverse mean curvature flow and the
  {R}iemannian {P}enrose inequality}, J. Differential Geom. \textbf{59} (2001),
  no.~3, 353--437.

\bibitem{KW}
J.~L. Kazdan and F.~W. Warner, \emph{Prescribing curvatures}, Differential
  geometry ({P}roc. {S}ympos. {P}ure {M}ath., {V}ol. {XXVII}, {S}tanford
  {U}niv., {S}tanford, {C}alif., 1973), {P}art 2, Amer. Math. Soc., Providence,
  R.I., 1975, pp.~309--319.

\bibitem{Lich}
A.~Lichnerowicz, \emph{Spineurs harmoniques}, C. R. Acad. Sci. Paris
  \textbf{257} (1963), 7--9.

\bibitem{Lohk}
J.~Lohkamp, \emph{Scalar curvature and hammocks}, Math. Ann. \textbf{313}
  (1999), no.~3, 385--407.

\bibitem{Minoo}
M.~Min-Oo, \emph{Scalar curvature rigidity of asymptotically hyperbolic spin
  manifolds}, Math. Ann. \textbf{285} (1989), no.~4, 527--539.

\bibitem{S}
R.~Schoen, \emph{Variational theory for the total scalar curvature functional
  for {R}iemannian metrics and related topics}, Topics in calculus of
  variations ({M}ontecatini {T}erme, 1987), Lecture Notes in Math., vol. 1365,
  Springer, Berlin, 1989, pp.~120--154.

\bibitem{SYNAS}
R.~Schoen and S.-T. Yau, \emph{Incompressible minimal surfaces,
  three-dimensional manifolds with nonnegative scalar curvature, and the
  positive mass conjecture in general relativity}, Proc. Nat. Acad. Sci. U.S.A.
  \textbf{75} (1978), no.~6, 2567.

\bibitem{SYsc}
\bysame, \emph{Existence of incompressible minimal surfaces and the topology of
  three-dimensional manifolds with nonnegative scalar curvature}, Ann. of Math.
  (2) \textbf{110} (1979), no.~1, 127--142.

\bibitem{SY1}
\bysame, \emph{On the proof of the positive mass conjecture in general
  relativity}, Comm. Math. Phys. \textbf{65} (1979), no.~1, 45--76.

\bibitem{SYsc2}
\bysame, \emph{On the structure of manifolds with positive scalar curvature},
  Manuscripta Math. \textbf{28} (1979), no.~1-3, 159--183.

\bibitem{SY2}
\bysame, \emph{Proof of the positive mass theorem. {II}}, Comm. Math. Phys.
  \textbf{79} (1981), no.~2, 231--260.

\bibitem{Wang}
X.~Wang, \emph{The mass of asymptotically hyperbolic manifolds}, J.
  Differential Geom. \textbf{57} (2001), no.~2, 273--299.

\bibitem{WY}
E.~Witten and S.-T. Yau, \emph{Connectedness of the boundary in the
  {A}d{S}/{CFT} correspondence}, Adv. Theor. Math. Phys. \textbf{3} (1999),
  no.~6, 1635--1655 (2000),.

\bibitem{Yau}
S.-T. Yau, \emph{Geometry of three manifolds and existence of black hole due to
  boundary effect}, Adv. Theor. Math. Phys. \textbf{5} (2001), no.~4, 755--767.




\end{thebibliography}
\end{document}